\input amstex
  \documentstyle{amsppt}
  \headline={\hss\tenrm\hss\folio}
  \vsize=9.0truein
  \hsize=6.5truein
  \parindent=20pt
  \baselineskip=22pt
  \NoBlackBoxes

   \topmatter
   %TITLE
	\title An Orbifold Relative Index Theorem \endtitle 

	 \author Carla Farsi, \
	  Department of Mathematics, \ 
	   University of Colorado, \
		395 UCB, \
		 Boulder, CO 80309--0395, USA. \ 
		  e-mail: farsi\@euclid.colorado.edu 
		   \endauthor
			   % ABSTRACT
				\abstract 
				 In this paper we prove a relative index theorem for
				 pairs of generalized Dirac operators on 
				 orbifolds which are the same at infinity. This
				 generalizes to orbifolds a celebrated theorem  of
				 Gromov and Lawson.
				  \endabstract
				   \endtopmatter

					\noindent {\bf 0.  Introduction.}
					 
					 Orbifolds, generalized manifolds that are locally
					 the quotient of an euclidean  space modulo a finite
					 group of isometries, were first introduced by Satake. In
					 the late seventies, Kawasaki proved an orbifold
					 signature formula, together with more general
					 index formulas, see [Kw1], [Kw2], [Kw3].  In [Fa1] we proved a $K$--theoretical index theorem
					 for orbifolds with operator algebraic means,
					  and  in [Fa2] and  [Fa3] we studied compact
					 orbifold spectral theory and defined orbifold eta
					 invariants. Other orbifold index formulas were proved in [Du], [V]. In [Ch] Chiang studied compact
					 orbifold heat kernels and harmonic maps, while in
					 [Stan], Stanhope established some interesting
					 geometrical applications of orbifold spectral
					 theory. 

					 In [Fa4] we continued our orbifold spectral
					 analysis started with [Fa2] and [Fa3]. In
					 particular we showed that on a non--compact
		complete almost complex $Spin^c$ orbifold which is sufficiently regular at infinity
(see Definition 2.1), generalized Dirac operators
					 are closed. This generalizes to orbifolds theorems
					 of Gaffney [Gn1], Yau [Y], and Wolf [W] for the manifold case.  We also proved an orbifold Divergence/Stokes Theorem.

					 Here we will use the results we proved
					 in [Fa4] to establish an orbifold
Gromov-Lawson, [GL], relative index theorem for non--compact
complete almost complex $Spin^c$  orbifolds
        which are sufficiently good at infinity
					 (see the definition right before the statement of Theorem 3.5). Further generalizations to
					 orbifolds of theorems in the spirit of the main
					 result in [An1], [An2] will be considered in a later
					 paper, [Fa5], together with an analytic proof
					 of the cobordism invariance of the index. We will also plan on investigating the question of whether every open complete orbifold is sufficiently good at infinity, which we suspect has a negative answer.

The relative index theorem applies to pairs of 
					 $+$  generalized Dirac operators which agree at infinity.
The topological index of such a pair is defined to be the difference of the index of two natural extensions of the given operators to closed manifolds, see Definition 4.1. Kawasaki's theorems for closed orbifolds makes possible to express
the pair topological index via local traces and orbifold characteristic classes. The
topological index is calculated by using pairs of parametrices that agree at infinity, c.f. Theorem 4.4.
					 
The analytical index of a $+$  generalized Dirac
					 operator, proved to be finite in Section 3, is
					 defined, classically, to be the dimension of the
					 kernel minus the dimension of the cokernel.
					 The analytical index of a pair is defined to be the difference of the analytical indices of the pair's constituents, see Definition 3.6. 

Our main result is Theorem 6.2, which asserts the equality of the analytical and topological indices of a pair of generalized Dirac operators that agree at infinity,
which generalizes to orbifolds the main theorem of [GL]. 
\proclaim{Theorem 6.2} Let $X$ be an even--dimensional non--compact 
complete almost complex Hermitian $Spin^c$ 
orbifold which is sufficiently good at infinity. Assume that a Hermitian connection is chosen on the dual of its canonical line bundle.
Let $D_i$ and $D_i^\pm$ be the 
generalized Dirac operators on $X$ with coefficients in the
proper Hermitian orbibundle  $E_i$ (with connection $\nabla^{E_i}$),  $i=1,2$. Suppose that $D_1 = D_2$ in a neighborhood $\Omega$ of infinity.
Assume that there exists a constant $k_0>0$ such that
$$
{\Cal R} \geq \, k_0 \, Id\, \hbox{ on } \Omega,
$$   
where ${\Cal R}$ is given as in Proposition 2.8. Then the topological index (as
in Definition 4.1) and the analytical index (as in Definition
3.6) of the pair $( {D}_1^+, { D}_2^+)$ coincide, that is
$$
ind_t ( {D}_1^+, { D}_2^+) = ind_a ( {D}_1^+, { D}_2^+) $$
\endproclaim

					 Our proof is a generalization to orbifolds of the
					 argument given for manifolds by Gromov and Lawson in [GL].
					 New techniques are mainly introduced to deal with
					 orbifold distance functions.

					 More in detail, the contents of this paper are as
					 follows.  In Section 1
					 we recall the definition of orbifold, orbibundles, and introduce orbifold generalized Dirac operators. In
					 Section 2 we recall properties of  generalized Dirac operators
					 on non--compact orbifolds from [Fa4]. In
					 Section 3 we prove that generalized Dirac operators (on
					 complete sufficiently regular orbifolds) which are
					 positive at infinity, have finite dimensional
					 kernels and cokernels. By using these results we can thus define the analytical index of a pair of generalized Dirac operators that agree at infinity as the difference of their analytical indices, see Definition 3.6. In
					 Section 4 we define the pair topological index by using a gluing technique and characterize it by using Kawasaki's index theorem for closed orbifolds, see Definition 4.1 and Theorem 4.2. In Theorem 4.4 we show how to compute the 
pair topological index by using parametrices that agree at infinity.   Section 5 is devoted to detailing some local properties of traces of generalized
Dirac
operators. Finally, in  Section 6 we state and prove our main result, the orbifold relative
index theorem.  

 In the sequel, all orbifolds and manifolds are
assumed to be even dimensional,
smooth, Hermitian, $Spin^c$, connected, and almost complex 
unless otherwise specified. All vector and orbibundles are assumed to be smooth and proper. We also assume that all of our orbifolds/manifolds are endowed with a fixed connection on the dual of their canonical line bundle $K^*$.
This allows us to define a \lq canonical\rq $\, Spin^c \,$ 
Dirac operator and, given an Hermitian orbibundle $E$ endowed with a connection, the \lq canonical\rq $\, Spin^c\, $ 
Dirac operator with coefficients in $E$. Both of these operators depend, in the $Spin^c$ case, on the choice of the selected connections, see [Du; Chapter 14], and [LM; Appendix D]. For the $Spin$ or complex case, the choice of the connection on $K^*$ is canonical.

					 I would like to thank the sabbatical program of the University of Colorado/Boulder, and the Mathematics Department of the University of Florence,
					 Italy, for their warm hospitality during the period
					 this paper was written. I would also like to
					 thank the  referee of pointing out the incorrect formulations of some of our results, and for his/her suggestions.  

					  \vskip 1em
					   \vskip 1em
						  \noindent {\bf 1.  Orbifolds, Orbibundles and
						  Dirac Operators.}
						   \vskip 1em
							
							 In this section we will review some
							 definitions and results that we will use throughout 
							 this paper. For generalities on orbifolds and operators on orbifolds, 
see [Kw1], [Kw2], [Kw3], [Ch], [Du]. 

							 An orbifold is a Hausdorff second countable
							 topological space $X$ together with an
							 atlas of charts 
							 $\Cal U  =\{ {(\tilde U_i}, G_i )| i \in {
							 I\} }$, with ${\tilde U}_i/G_i = U_i$ open  and
							with projection $\pi_i: {\tilde U}_i \to U_i$, $i\in I$, 
							 satisfying the following properties 
							 \roster
							 \item  If  two charts $U_1$ and $U_2$
							 associated to pairs $({\tilde U}_1, G_1)$,
							 $({\tilde U}_2, G_2)$ of $\Cal U$,  are such that $U_1\subseteq U_2$, then there exists a smooth
							 open embedding $\lambda$:  ${\tilde U}_1\to {\tilde U}_2$ and a homomorphism
							 $\mu$: $G_1\to G_{2}$ such that $\pi_1=
							 \pi_{2}\circ \lambda$ and $\lambda \circ
							 \gamma =\mu(\gamma)\circ \lambda$, $\forall
							 \gamma \in G_1$.
							 \item The collection of the open charts
							 $U_i$, $i\in I$, belonging to the atlas
							 $\Cal U$ forms a basis for the topology on
							 $X$. 
							 \endroster
							 We will call an orbifold atlas as above a
							 standard orbifold atlas.

							 For any $x$ point of $X$, the isotropy
							 $G_x$ of $x$ is well defined, up to
							 conjugacy, by using any local coordinate
							 chart. The set of all points $x\in X$ with
							 non--trivial stabilizer, $\Sigma(X)$, is
							 called the singular locus of $X$. $\Sigma(X)$ is a set of codimension at least 2, see e.g.
							 [Ch].  Note that $X-\Sigma(X)$ is a 
 manifold.
								
								If we now endow $X$ with a countable
								locally finite orbifold atlas $\Cal F$, $\Cal F  =\{(
								{\tilde U_i}, G_i )| i \in {\Bbb N\} }$,
							then 	by standard theory
								there exists a smooth partition of unity
								$\eta=\{ \eta_i \}_{i\in \Bbb N}$
								subordinated to $\Cal F$, [Ch]. This in
								particular means that, for any $i \in
								\Bbb N$, $\eta_i$ is a smooth function
								on $U_i$ (i.e., its lift to any chart of
								a standard orbifold atlas is smooth),
								the support of $\eta_i$ is included in
								an open subset  $U_i'$ of $U_i$, and  $\cup U_i'=X$. We will call any $\eta$ as above an $\Cal F$-partition of unity. 

								Let $E$ be an Hermitian orbibundle (with connection $\nabla^E$)
								over the orbifold $X$. (For the precise
								definition see [Kw1], [Kw2], [Kw3],
								[Ch].) In particular $E$ is an orbifold
								in its own right; on an orbifold chart
								$U_1$ associated to a pair $({\tilde
								U}_1, G_1)$ of a standard orbifold atlas
								$\Cal U  =\{ ({\tilde U_i}, G_i )| i \in {
								I\} }$ of $X$, $E$ lifts to a
								$G_1$-equivariant bundle. Standard
								orbifold atlases on $X$ can be used
								to provide standard orbifold atlases on
								$E$. 

								If $E$ is an orbibundle
								over the orbifold $X$, a section $s:X\to
								E$ is called a smooth orbifold section
								if for each chart $U_i$ associated to a
								pair $({\tilde U}_i, G_i)$ of a standard
								orbifold atlas $\Cal U  =\{ ({\tilde
								U_i}, G_i) | i \in { I\} }$ of $X$, we
								have that $s|_{U_i}: U_i \to E|_{U_i}$
								is covered by a smooth $G_i$--invariant
								section ${\tilde s}|_{{\tilde U}_i}:
								{\tilde U}_i \to {\tilde E|_{{\tilde U}_i}}$.
								Given an orbibundle $E$ over $X$,
								we will denote by 
								${\Cal C}^\infty(X,E)$ the space of all
								smooth sections of $E$, and by
								${\Cal C}^\infty_c(X,E)$  the space of
								all smooth sections with compact
								support.
								Classical orbibundles over $X$
								are the tangent bundle $TX$, and the
								cotangent bundle $T^*X$ of $X$. We can
								form orbibundle tensor products by
								taking the tensor products of their local
								expressions in the charts of a standard
								orbifold atlas.  

								Define an inner product between sections                 of ${\Cal C}^\infty (X,E)$  (or
${\Cal C}^\infty_c(X,E)$) by a the
								following formula (c.f., [Ch; 2.2a]
								$$
								(\sigma_1, \sigma_2)=
								\sum_{i=1}^{+\infty} \frac1{|G_i|}
								\int_{\tilde U_i}
								 {\tilde \eta_i} (\tilde x_i ) <{\tilde
								 \sigma_1}(\tilde x_i ), {\tilde
								 \sigma_2}(\tilde x_i )> dv(\tilde x_i
								 ),
								 $$ 
								 where $\eta=\{ \eta_i \}_{i\in \Bbb
								 N}$, is a $\Cal F$-partition of unity
								 subordinated to the locally finite
								 orbifold cover $\Cal F  =\{ ({\tilde
								 U_i}, G_i) | i \in {\Bbb N} \}$, and
								 $<,>$ is a $G_i$--invariant product on $\tilde
								 E$. (Note that, with slight abuse of
								 notation, we used $\,\tilde {}\,$ to denote
								 lift to ${\tilde U_i}$.)

								 We will now review the construction of the
								 generalized Dirac operator
with coefficients in an Hermitian orbibundle 
$E$ (with connection $\nabla^E$)  over a (compact or not) orbifold
								 $X$, see [Du; Sections 5 and 12], [Kw2], [BGV], and [LM; Appendix D] in the manifold case.  First of all, $X$  admits a
								 $Spin^c$-principal tangent orbibundle,
								 $Spin^c(TX)$, with, in our hypotheses, a  canonical $Spin^c$
								 orbifold connection $\nabla^c$. 
								 Let $\Delta^{\pm,
								 c}$ be the half  $Spin^c$
								 representations (recall that the $X$ is even dimensional).  Then we have two
								 orbibundles
								 $$
								 \Delta^{\pm, c} (TX) = Spin^c(TX)
								 \times_{ Spin^c} \Delta^{\pm, c}, $$
								 with induced connections $\nabla^{\pm,
								 c}$, from $\nabla^c$; $\nabla^{\pm,c} : {\Cal C}^\infty_c (X, \Delta^{\mp, c} (TX))
\to {\Cal C}^\infty_c (X, T^*X \otimes \Delta^{\mp, c} (TX))$.

								 The Clifford module structure on
								 $\Delta^{\pm, c}$ defines Clifford
								 multiplications
								 $$
								 m_{\pm}: TX \otimes_{\Bbb R}
								 \Delta^{\pm, c} (TX) \to \Delta^{\mp,
								 c} (TX)
								 $$ 
								 Then the
								 generalized $\pm $ Dirac operator with
								 coefficients in $E$, $d_E^{\pm, c}$, 
								 $$
								 d_E^{\pm, c}: {\Cal C}^\infty_c (X,
								 \Delta^{\pm, c} (TX)\otimes_{\Bbb C} E)
								 \to
								 {\Cal C}^\infty_c (X, \Delta^{\mp, c}
								 (TX)\otimes_{\Bbb C} E)
								 $$
								 is defined by
								 $$
								 d_E^{\pm, c}= M \circ \left( \nabla^{\pm, c}\otimes Id + Id \otimes \nabla^{E}\right),
								 $$
								 where $M$ denotes the map induced by Clifford multiplication and TX has been identified with $T^*X$ via
								 the orbifold metric. 
We will also use the notation  $\Cal S$ for
								 the orbifold $Spin^c$ bundle
								 $(\Delta^{+, c}\oplus \Delta^{-, c}) (TX)$, and ${\Cal S}\otimes E$ or $\Cal E$  for $(\Delta^{+,
								 c}\oplus \Delta^{-, c})
								 (TX)\otimes_{\Bbb C} E$. 
We will define $D_E$, the generalized
								 Dirac operator on $X$ with coefficient
								 in $E$, to be $(d_E^{+, c}+ d_E^{-, c})$.

								 \vskip 1em
								 \noindent {\bf 2.  Dirac Operators on Non--Compact Complete Orbifolds.} \vskip 1em

								 On an orbifold $X$ satisfying  our hypotheses (but not necessarily
								 compact), the generalized Dirac
 operator $D_E$  
with coefficients in the
								 orbibundle $E$, as
	 defined in the Section 1, is given by 
								 $$
								 D_E: {\Cal C}^\infty_c (X, (\Delta^{+,
								 c}\oplus \Delta^{-, c})
								 (TX)\otimes_{\Bbb C} E) \to {\Cal
								 C}^\infty_c (X, (\Delta^{-, c}\oplus
								 \Delta^{+, c}) (TX)\otimes_{\Bbb C} E)
								 $$
								 $$
								 D_E= M \circ \left( (\nabla^{+, c} +
								 \nabla^{-, c})\otimes Id + Id \otimes
								 \nabla^{E} \right).
								 $$

								 On orbifold charts,								 $D_E$ has the following local
								 expression ${\tilde D}_E$. Let $\Cal U
								 =\{({\tilde U_i}, G_i) | i \in { I\} }$,
								 with ${\tilde U}_i/G_i = U_i$ be a
								 standard orbifold atlas. On a local
								 chart ${\tilde U}_i$, $i\in I$ fixed, we have
								 $$
								 \Delta^{\pm, c} (T{\tilde U_i}) =
								 Spin^c(T{\tilde U_i}) \times_{ Spin^c}
								 \Delta^{\pm, c}, $$
								 with induced $G_i$--invariant
								 connections $\nabla^{\pm, c}$, from
								 $\nabla^c$.
								 The Clifford module structure on
								 $\Delta^{\pm, c}$ defines Clifford
								 multiplications
								 $$
								 m_{\pm}: T{\tilde U_i} \otimes_{\Bbb R}
								 \Delta^{\pm, c} (T{\tilde U_i}) \to
								 \Delta^{\mp, c} (T{\tilde U_i}).
								 $$ 
								 On $\tilde E$, the lift of $E$,  we have the 
								 $G_i$--invariant connection
								 $\nabla^{\tilde E}$. 
								 Then the generalized $\pm $
								 Dirac operators with coefficients in
								 $E$, ${\tilde d}_E^{\pm, c}$, 
								 $$
								 {\tilde d}_E^{\pm, c}: {\Cal
								 C}^\infty_c ({\tilde U_i}, \Delta^{\pm,
								 c} (T{\tilde U_i})\otimes_{\Bbb C}
								 {\tilde E}) \to
								 {\Cal C}^\infty_c ({\tilde U_i},
								 \Delta^{\mp, c} (T{\tilde
								 U_i})\otimes_{\Bbb C} E)
								 $$
								 is given by
								 $$
								 {\tilde d}_E^{\pm, c}= M \circ \left(
								 \nabla^{\pm, c}\otimes Id + Id
								 \otimes \nabla^
								 {\tilde E} \right),
								 $$
								 where $M$ is induced by Clifford multiplication and $T{\tilde U_i}$ has been identified
								 with $T^*{\tilde U_i}$ via the
								 $G_i$--invariant metric. 
								 Also, ${\tilde D}_E$, the generalized
								 Dirac operator on $X$ with coefficient
								 in $E$, is given by ${\tilde d}_E^{+,
								 c}+ {\tilde d}_E^{-, c}$ on ${\tilde U}_i$. 

								 If $e_1, \dots, e_n$ is an orthonormal
								 local basis for the space $T{\tilde
								 U_i}$ at a point $\tilde x$, then ${\tilde
								 D}_E$
								 has local expression
								 $$
								 {\tilde D}^{ E} =
								 \sum_{k=1}^{n} e_k {\tilde
								 \nabla}^{E}_{e_k},
								 $$
								 where 
								 $$
								 {\tilde \nabla}^{E}= ({\nabla}^{+, c}+
								 {\nabla}^{-, c} ) \otimes 1 + 1 \otimes
								 {\nabla}^{\tilde E}.
								 $$

								 In analogy with the manifold case,
								 see [GL], [W], [Gn1], [LM], [Y], one can show that $D_E$ is symmetric,
								 whenever $X$          
								 sufficiently regular at
								 infinity.

								 \proclaim{Definition 2.1 } Let $X$ be a
								 non--compact complete orbifold. Then we say that
								 $X$ is sufficiently regular at infinity
								 if, for any neighborhood  $\Omega
								 \subseteq X$ of infinity, there exists
								 a compact domain $K_\Omega$ with boundary strictly included in $\Omega$,
								 on which the Divergence and Stokes' Theorems hold.  
								 \endproclaim

								 For a compact orbifold
								 without boundary, the Divergence
								 Theorem holds,
								 [Ch]. See also [C] for other
								 results. Sufficient regularity also holds in
								 the case of a product end, by an
								 adaptation of Chiang's method, see
								 [Ch]. In general, ours seems to be a
								 very reasonable assumption to make,
								 which will be certainly satisfied in
								 many cases of interest.  For Sobolev inequalities of Gallot type involving domains, see [N].

								 \proclaim{Theorem 2.2 [Fa4; Theorem
								 2.2]} Let $X$ be a
								 non--compact complete orbifold which is
								 sufficiently regular at infinity, and
								 let $E$ be a Hermitian orbibundle over $X$.
								 Let $D_E$ be the generalized
								 Dirac operator with coefficients in
								 $E$, as defined above. 
Then $D_E$ is symmetric, i.e.,    
								 $$
								 (D_E \sigma_1, \sigma_2)= (\sigma_1,
								 D_E\sigma_2), \quad \forall \sigma_1, \sigma_2
\in {\Cal C}^\infty_c (X,
{\Cal S}\otimes_{\Bbb C} E),
$$
 where $(,)$ denotes the inner product
					defined in the introduction. (Actually, to obtain the above equality, it is enough to assume that only one of the sections $\sigma_1, \, i=1,2$ has compact support.)
\endproclaim

								   Now complete the space ${\Cal
								   C}^\infty_c (X, {\Cal E})$, ${\Cal
								   E}= {\Cal S}\otimes_{\Bbb C} E$, ${\Cal S}$
								   $Spin^c$--bundle on $X$, $E$ Hermitian   orbibundle (with connection $\nabla^E$) over $X$, with respect to
								   the norm 
								   $$
								   \Vert \sigma \Vert_X= \sqrt{ <\sigma,
								   \sigma>}=
								   \left(\sum_{i=1}^{+\infty}
								   \frac1{|G_i|} \int_{\tilde U_i}
								   {\tilde \eta_i} (\tilde x_i )
								   <{\tilde \sigma_1}(\tilde x_i ),
								   {\tilde \sigma_2}(\tilde x_i )>
								   dv(\tilde x_i )
\right)^{\frac12}. 
								   $$
								  We thus obtain the
								   ${\Cal L}^2$--space ${\Cal L}^2(X,
								   {\Cal E})$.  The generalized Dirac operator
								   $$
								   D_E: {\Cal C}^\infty_c (X, {\Cal E})
								   \to {\Cal C}^\infty_c (X, {\Cal E}) 
								   $$
								   has two natural extensions, min and max, see [Fa4],  as an
								   unbounded operator
								   $$
								   D_E: {\Cal L}^2 (X, {\Cal E}) \to
								   {\Cal L}^2 (X, {\Cal E}). 
								   $$

								   \proclaim{Theorem 2.3 [Fa4; Theorem 3.1] } Let $X$ be a non—-compact
complete orbifold which is
sufficiently regular at infinity, and
let $E$ be a Hermitian orbibundle (with connection $\nabla^E$)  over
								   $X$, and  let $D_E$ be the generalized Dirac
								   operator with coefficients in $E$.
								   Let ${\Cal D}( D_E^{MIN})$ be the domain of the min extension of $D_E$,
								   and  $ {\Cal D}( D_E^{MAX})  $ be the domain of the max extension of $D_E$. Then
								   $$
								   {\Cal D}( D_E^{MIN})= {\Cal D}(
								   D_E^{MAX}).  
								   $$
								   \endproclaim
The following very useful proposition was proved in [Fa4]; it
is a generalization of results proved by Gaffney and Yau for manifolds,
see [Gn2], [Y].
\proclaim{Proposition 2.4 [Fa4; Proposition 3.4]} Let $X$
									 be a non--compact complete orbifold
									 which is sufficiently regular at
									 infinity, and let $y_0 \in
									 X-\Sigma(X)$ be a fixed 
									 point of $X$. Then there exists a
									 sequence of continuous functions
									 $b_k$, $k\in \Bbb N$, with
									 \roster
									 \item $ b_k: X\to [0,1]$
									 \item $b_k=1$ on $B_k= \{ y\in X|
									 \rho(y) =d(y,y_0)\leq k\}$. 
									 \item The support of $b_k$ is
									 contained in ${\overline B}_{2k}$.
									 \item The function $b_k$ is
									 differentiable almost everywhere
									 and at points of differentiability
									 we have
									 $$
									 \Vert \nabla (b_k) \Vert^2 \leq
									 \frac{M^2}{k^2}, \quad k \in \Bbb N 
									 $$
									 \endroster
									 \endproclaim

In [Fa4] we also proved the following results.

\proclaim{Theorem 2.5 [Fa4; Theorem 4.1]} Let $X$ be  non--comapct
complete orbifold which  is sufficiently
  regular at infinity, and let $E$  be a Hermitian orbibundle (with connection $\nabla^E$)  over $X$.
									   Let
									   $$
				D_E: {\Cal C}^\infty_c (X,
									   \Cal E ) \to {\Cal
				 C}^\infty_c (X, \Cal E),
									   $$
									   be the generalized Dirac operator
						   on $X$ with coefficients in $E$.
									   Then
   $D_E(\sigma)=0$ if and only if
									   $D_E^2(\sigma)=0$ for any $\sigma
									   \in {\Cal D}(D_E)$.\endproclaim

\proclaim{Proposition 2.6 [Fa4; Proposition 6.1]}
											  Let $X$ be a non--compact
											  complete $Spin^c$ orbifold
										which is 	  sufficiently regular at
											  infinity, and let $\Cal S$ be  the $Spin^c$ bundle of $X$. 
											  Then for any two
											  sections $\sigma_j$,
											  $j=1,2$ in ${\Cal
											  C}^\infty (X, \Cal S)$, at
											  least one of which with
											  compact support, we have
											  $$
											  \int_X <\Delta \sigma_1,
											  \sigma_2 > dv =\int_X <
											  \nabla\sigma_1,
											  \nabla\sigma_2 > dv 
											  $$ 
											  \endproclaim

\proclaim{Proposition 2.7 [Fa4; Proposition 6.3]}
 Let $X$ be an orbifold.  
 Let $D$ be
the Dirac
											  operator on $X$, and let $\Delta$
											  be the $Spin^c$ Laplacian.
											  Then
											  $$
	 D^2= \Delta+ {\Cal R}
											  $$
  where ${\Cal R}$ is given below  (c.f. [Du; Theorem 6.1], [LM; Theorem D12]
for the manifold case),
$$
\Cal R= \frac14\,  k\,  + \frac12 c(K^*),
$$
where $k$ is the scalar curvature,  and $ c(K^*)$ denotes the
Clifford multiplication of the curvature 2—form of the fixed  
connection on the line bundle $K^*$. 
\endproclaim
When  $D_E$ is the generalized Dirac operators on $X$ with coefficients in the
proper Hermitian orbibundle (with connection $\nabla^E$)  $E$, then the formulas above 
become
$$ D_E^2= \Delta_E+ {\Cal R_E}, 
$$
$$
\Cal R_E= \frac14\,  k\,  + \frac12 c(K^*) + c(E),
$$
where $c(E)$ is the Clifford multiplication of the curvature 2—form of the fixed  
connection $\nabla^E$ on $E$.   
When $X$ is $Spin$, we can assume that $ c(K^*) = c(E)=0$, and ${\Cal R}_E = \Delta_E + \frac 14 k$. 

\proclaim{Theorem 2.8 [Fa4; Theorem 6.4]} Let
$X$ be a complete,
non--compact, $Spin^c$
											  orbifold which  is
											  sufficiently regular at
											  infinity. If $D$ is the
											  Dirac operator on $X$ with
											  coefficients in the
											  $Spin^c$ bundle $\Cal S$,
											  then
											  the domain $\Cal D$ of the
											  unique self-adjoint
											  extension of $D$ is
											  exactly
											  $$
											  {\Cal L}^{1,2}(X, {\Cal
											  S}), \quad \hbox{that is,}
											  $$
											  the completion of ${\Cal
											  C}^\infty_c (X, {\Cal S})
											  $  in the norm 
											  $$
											  \Vert \sigma\Vert_1^2=
											  \int_X\left( <\sigma,
											  \sigma>+ <\nabla\sigma,
											  \nabla\sigma>\right)dv =
											  \int_X\left( <\sigma,
											  \sigma>+ <\Delta_\sigma,
											  \sigma>
											  \right)dv
											  $$
											  Moreover, for every
											  $\sigma\in \Cal D$,
											  $$
											  \Vert D\sigma\Vert^2_X =
											  \Vert \nabla\sigma \Vert^2_X
											  + ({\Cal R} \sigma,
											  \sigma), 
											  $$
											  where $\Vert\, \Vert_Y$
denotes the ${\Cal L}^2(Y,
											  \Cal E)$, $Y\subseteq
											  X$,
											  norm, $\Cal R$ is is given as in Proposition 2.7, and $(\, , )$ the
											  ${\Cal L}^2$ inner
											  product.
											  \endproclaim
							   \vskip 1em
								\noindent {\bf 3. Dirac Operators and
								Green's Operators.}
								 \vskip 1em

								 As pointed out by Agmon in [A; Section
								 6], Friedrichs' Lemma is a local
								 result. Hence we can also assume it holds for orbifolds. (This
								 follows from the local version of the manifold lemma applied
								 to covers of orbifold charts of a locally finite
								 orbifold cover.) 

								 \proclaim{Theorem 3.1 (Friedrichs' Lemma for Orbifolds)} Let
								 $X$ be an
orbifold. Let $\Cal
								 S$ be the $Spin^c$ bundle of $X$, and
								 let  $D$ be the Dirac operator on $X$, $D: {\Cal C}^\infty_c (X, {\Cal S}) \to
{\Cal C}^\infty_c (X, {\Cal S}). $
Let $\Omega$ be an open set in $X$ and
let $K$ be a compact subset of
$\Omega$. Let $k\in \Bbb N$. Then there
exists a constant $C$, depending only on
$K$, $\Omega$, and $k$, such that, for
								 any $\sigma \in   {\Cal C}^\infty_c (X,
								 {\Cal S}|_\Omega)$ with $D\sigma =0$,
								 we have
								 $$
								 \Vert \sigma\Vert_{{\Cal C}^k, K} \leq
								 \, C\, \Vert  \sigma\Vert_{{\Cal L}^2(\Omega,
								 \Cal S)},
								 $$ 
								 where $\Vert \, \Vert_{{\Cal C}^k, K}$
								 is the uniform ${\Cal C}^k$  norm for sections on $K$.
\endproclaim

								 \noindent {\it Proof.} See 
[GL; Theorem 3.7] for the
								 manifold case. Endow $X$ with a
								 countable locally finite orbifold atlas
								 $\Cal F$, $\Cal F  =\{( {\tilde U_i},
								 G_i )| i \in {\Bbb N\} }$,  with
								 associated  partition of unity $\eta=\{
								 \eta_i \}_{i\in \Bbb N}$ subordinated
								 to $\Cal F$, [Ch].  Suppose that $U_i
								 \cap K \not= \emptyset $ for only
								 $i=1,\dots, \ell$. We can also choose
								 our atlas so that $\cup_{i=1}^\ell U_i
								 \subseteq \Omega$.
								 Then by the manifold local version of
						Friedricks' Lemma [A; Section 6],
modulo multiplication by $|G_i|$, $i=1,\dots, \ell$, we obtain the claim. $\qed$ 

								 We can now prove that Dirac operators which are positive at
								 infinity, have finite dimensional
								 kernels and cokernels. 

								 \proclaim{Theorem 3.2} Let $X$ be a
								 non--compact complete  orbifold
which is sufficiently regular at
infinity. Let $\Cal S$ be the $Spin^c$
 bundle of $X$, and let  $D$ be the
								 Dirac operator on $X$, 
$D: {\Cal C}^\infty_c (X, {\Cal S}) \to
					 {\Cal C}^\infty_c (X, {\Cal S}). $
								 Assume that there exists a compact subset $K$ of $X$ such that
								 $$
								 {\Cal R} \geq k_0\, Id 
								 $$
								 in $X-K$, where $\Cal R$ is 
given as in Proposition 2.7.
								 Then there exists an integer $d$
								 depending only on $D$ on a neighborhood
								 of
								 $K$ and on $k_0$, such that
								 $$
								 \hbox {dim (Ker} D) \leq d.
								 $$
								 In particular, if the dimension of $X$
								 is even, and $D=D^+ \oplus D^-$ 
 (see Section 1), then
								 $$
								 \hbox {dim (Ker} D^+) +\hbox {dim (Ker}
								 D^-) \leq d.
								 $$
								 \endproclaim

								 \noindent {\it Proof.} Since $K$ is
								 compact, there exists a $k_1>0$ such
								 that $k_1 Id \geq -\Cal R$ in $K$. Now
								 let $\sigma \in {\Cal L}^2(X, \Cal S)$,
								 with
								 $D\sigma =0$. Firstly, $\sigma \in
								 {\Cal C}^\infty(X, \Cal S)$, since
								 since solutions of elliptic equations on
  manifolds are smooth by a local argument.  Then by Propositions 2.6								 and 2.7, we have, 
 for $\sigma$ with compact support, 
								 $$
								 \Vert \nabla \sigma \Vert^2_X
								  + ({\Cal R}\sigma, \sigma) = 0,
								 $$
								 where $\Vert \, \Vert_Y$, $Y\subseteq
								 X$, is the
								 ${\Cal L}^2(Y,\Cal S)$ norm.
								 Since ${\Cal R} \geq k_0 $ on $X-K$, we
								 have
								 $$
								 \Vert \nabla \sigma \Vert^2_X
								  + \int_K <{\Cal R} \sigma,
								 \sigma> + k_0 \int_{X-K}  \Vert \sigma
								 \Vert^2 \leq 0.
								 $$
								 Therefore, 
								 $$
								 \Vert \nabla \sigma \Vert^2_X + k_0
								 \Vert \sigma \Vert^2_{X-K} \leq \Vert
								 \sigma \Vert^2_{K}.
								 $$

								 We now make the assumption $\Vert
								 \sigma \Vert^2_X = \Vert \sigma
								 \Vert^2_K +
								 \Vert \sigma \Vert^2_{X-K}=1$.  By adding
								 $k_0 \Vert \sigma \Vert^2_K$ to both
								 sides, and dividing by $k_0 + k_1$, the
								 above inequality becomes
								 $$
								 \frac{\Vert \nabla \sigma \Vert^2_X
								 }{k_0+k_1} + \frac{k_0}{k_0+k_1} \leq 
								 {\Vert \sigma \Vert^2_K }. \tag 3.1
								 $$ 
								 Fix a neighborhood $\Omega$ of $K$
								 and let $C$ be the constant appearing
								 in Theorem 3.1 for $k=1$.  Fix
								 $\epsilon >0$ and an $\epsilon$ dense
								 subset
								 $\{ x_s\}_{s=1, \dots, d}$ of $K$ so
								 that every point in $K$ is within
								 distance
								 $\epsilon$ of some $x_s$. Let
								 $H=$Ker$(D)$ on ${\Cal L}^2(X, \Cal S)$
								 and suppose  by contradiction that dim$(H)>d$.  Then there
								 exists $\sigma \in H$ such that $\Vert
								 \sigma \Vert_X =1$, and $\sigma(x_s) =0$,
								 $s=1,\dots, q$. By Theorem 3.1,
								 applied to
								 $K$ and $\Omega$, and 
								 the Mean Value Theorem
								 applied on lifts of orbifold charts, we have
								 $$
								 \Vert \sigma (x)\Vert \leq sup_K \Vert
								 \nabla(\sigma)\Vert \, \epsilon \leq 
								 \, C \, \Vert \sigma\Vert_{\Omega}\, \epsilon.
								 $$ 
								 But this contradicts
								 $(3.1)$. (Note that we can assume that
								 $\sigma$ has compact support since we are
								 only interested in  behavior near
								 $K$.) $\qed$

								 By applying Glazman's variational lemma
								 below, see [Ku; Proposition 3.4] or [Sh], 
								 we can obtain information on the
								 spectrum of the Dirac operator.

								 \proclaim{Proposition 3.3} Let $A$
								 be a self-adjoint Hilbert operator
								 that is semi--bounded from below. Let
								 $N_h(\lambda)$ denote the number of
								 eigenvalues in $(-\infty, \lambda]$,
								 with multiplicity, and counting points
								 of the continuous spectrum as points
								 with infinite multiplicity. Then
								 $$
								 N_h(\lambda) =\hbox{ sup}_H \,
								 \hbox{dim}\, H ,
								 $$     
								 where the supremum is taken over all
								 the subspaces $H$ which are such that
								 $$
								 <Ah, h> \,\leq \, <h,h>, \, \forall h \in H
								 .
								 $$ 
								 \endproclaim

								 \proclaim{Corollary 3.4} Let $X$,
								 $D$, and $\Cal S$ be as in Theorem
								 3.2. Then $D^2$, and consequently
								 $D$, have essential spectrum separated from 0.
								 \endproclaim

								 We will now show that Dirac operators
								 admit Green's operators, c.f. [GL; Theorem 3.7], and [Ku; Proposition 3.3 and 3.4] for the manifold case. We assume here that the orbifold $X$ is
								 sufficiently good at infinity. This
								 means that, for any neighborhood of
								 infinity $\Omega$, we can chop off $X$
								 along an orbifold hypersurface,
								 $O_\Omega $, which is the boundary of a
								 neighborhood of infinity included in	 $\Omega$. We also assume that $\Omega$ is of product type near $O_\Omega$ 
\proclaim{Theorem 3.5} Let $X$ be a non--compact complete $Spin^c$ orbifold
								 which is sufficiently good at
								 infinity. Let $\Cal S$ be the $Spin^c$
								 bundle of $X$, and let  $D$ be the
								 Dirac operator on $X$,
$D: {\Cal C}^\infty_c (X, {\Cal S}) \to
{\Cal C}^\infty_c (X, {\Cal S}).  $
								 Assume that there exists a compact
								 subset $K$ such that
								 $$
								 {\Cal R} \geq k_0\, Id 
								 $$
								 on $X-K$, where $\Cal R$ is given as in Proposition 2.7. Assume that $H$ is the
								 finite--dimensional kernel of $D$ on
								 ${\Cal L}^2 (X, {\Cal S})$, and let
								 $H^\perp$ be its orthogonal complement.
								 Then there is a $\alpha>0$
								 such that 
								 $$
								 \Vert D\sigma \Vert^2_X \geq \alpha^2
								 \Vert \sigma \Vert^2_X, \, \forall
							\sigma \in	 H^\perp.
								 $$ 
								 Thus the operator $D$, and also, since								 $X$ is even-dimensional, the operator 
								 $D^\pm$, admits bounded Green's
								 operators. 
								 \endproclaim

								 \noindent{\it Proof.} (C.f. [GL; Proof of Theorem 3.7] and [Ku; 3.3 and 3.4].) First, we claim that $D$ has only point spectrum. This is proved by the same argument used in [Ku; 3.3.]. Note that Kucerovsky's proof works verbatim in the orbifold case too, since Rellich's lemma applies in the sufficiently good case (in fact, one can double the domain to get an orbifold, and then apply Rellich's Lemma for closed orbifolds proved in [Fa2]). Now, let $E_\lambda$
								 be the $\lambda$--eigenspace of
								 $D$ on ${\Cal L}^2 (X, {\Cal S})$ with
								 eigenvalue $\lambda$. It will suffice
								 to prove that there exists an
								 $\alpha>0$ such that the space 
								 $$
								 H_\alpha = \oplus_{|\lambda| \leq
								 \alpha} E_\alpha
								 $$
								 is finite-dimensional. We will now
								 proceed as in the proof of Theorem 3.2.
								 By Theorem 2.8, we have
								 $$
								 (D^2 \sigma, \sigma) -({\Cal R}
								 \sigma, \sigma) - \Vert \nabla \sigma
								 \Vert^2_X =0, \, \forall \sigma \in
								 {\Cal D}(D).
								 $$
								 Hence, for $\sigma \in H_\alpha$,
								 $D\sigma =\lambda\sigma$, with
								 $|\lambda| \leq \alpha$,  so we get
								 $$
								 \lambda^2 \Vert \sigma \Vert_X^2 - ({\Cal
								 R} \sigma, \sigma) - \Vert  \nabla
								 \sigma \Vert^2_X =0.
								 $$
								 Since for some $k_1\in \Bbb R$, $-{\Cal
								 R} \geq k_1 Id$ on $K$, we have that,
								 for any $\sigma \in H_\alpha$,
								 $$
								 \lambda^2 \Vert \sigma \Vert^2_X
								 +\int_K k_1 <\sigma, \sigma> dv \geq
								 \int_{X-K}
								 <{\Cal R}\sigma, \sigma> dv +
								  \Vert \nabla \sigma \Vert^2_X
								 $$
								 As $k_0\in \Bbb R$ is such that $k_0 \,Id
								 \leq {\Cal R}$ on $X-K$, we get
								 $$
								 \lambda^2 \Vert \sigma \Vert^2_X + k_1
								 \Vert \sigma \Vert^2_K\geq 
								 k_0 \Vert \sigma \Vert^2_{X-K} + 
								 \Vert \nabla \sigma \Vert^2_X, \,
								 \forall \sigma \in  H_\alpha. 
								 $$
								 Now replace  $\Vert \sigma
								 \Vert^2_{X-K} $ by $\Vert \sigma
								 \Vert^2_{X} -\Vert \sigma \Vert^2_{K}$
								 to get,  
								 $$
								 \lambda^2 \Vert \sigma \Vert^2_X + k_1
								 \Vert \sigma \Vert^2_K\geq 
								 k_0 (\Vert \sigma \Vert^2_{X} -\Vert
								 \sigma \Vert^2_{K})+  \Vert \nabla
								 \sigma \Vert^2_X, \, \forall \sigma \in
								 H_\alpha, 
								 $$
								 $$
								 \lambda^2 \Vert \sigma \Vert^2_X +
								 (k_0+ k_1) \Vert \sigma \Vert^2_K\geq 
								 k_0 \Vert \sigma \Vert^2_{X} + \Vert
								 \nabla \sigma \Vert^2_X, \,\forall
								 \sigma \in  \sigma \in H_\alpha. 
								 $$
								 Because $\sigma \in H_\alpha$,
								 $|\lambda| \leq \alpha$, and so
								 $\lambda^2 \leq \alpha^2$, which
								 implies
								 $$
								 (k_0 + k_1)  \Vert \sigma \Vert^2_K
								 \geq  (k_0- \alpha^2) \Vert \sigma
								 \Vert^2_X +  \Vert \nabla \sigma
								 \Vert^2_X, \,\forall \sigma \in
								 H_\alpha. 
								 $$
								 If we choose $\alpha>0$ such that
								 $\alpha^2 < < k_0$, we have
								 $$
								 \frac{k_0- \alpha^2}{k_0+ k_1} \leq
								 \frac{\Vert \sigma \Vert^2_K}{\Vert
								 \sigma \Vert^2_X}, \, \forall \sigma
								 \in  H_\alpha. \tag 3.2
								 $$
								 Now choose a parametrix $Q$ for $D$,
								 [Kw2], so that $QD=Id -T$, with $T$
								 smoothing. Let $\rho$ denote the
								 restriction to $K$, and let $\breve T =
								 \rho \circ T $. Then for any $ \sigma
								 \in H_\alpha$, we have 
								 $$
								 {\breve T} \sigma = \rho \sigma - \rho
								 QD \sigma.
								 $$    
								 Moreover, for any $\sigma \in
								 H_\alpha$, $\Vert D\sigma\Vert_X \leq
								 \alpha \, \Vert\sigma \Vert_X$. Set $q=\Vert Q
								 \Vert$; from (3.2) we get,
								 $$
								 \Vert {\breve T} \sigma \Vert_X \geq
								 \Vert \rho \sigma \Vert_X - \alpha \,  q \,  
								 \Vert \sigma \Vert_X \geq (\alpha_1
								 -\alpha q) \Vert \sigma \Vert_X, \forall
								 \sigma \in H_\alpha,
								 $$ 
								 where $\alpha_1 =
								 \frac{k_0-\alpha^2}{k_0 + k_1}$. Choosing $\alpha$ sufficiently small, $
								 \Vert {\breve T} \sigma \Vert_X \geq
								 \alpha_2 \Vert \sigma \Vert_X$,
								 $\forall \sigma \in H_\alpha$.
								 Therefore $H_\alpha$ is
								 finite-dimensional since ${\tilde T}$
								 is a compact operator. $\qed$

We will now define  the											 analytical											 index of a pair											 of generalized Dirac									
	 operators that		 agree at 											 infinity, thus							 generalizing a construction of								 Gromov and	Lawson.									 We start with a
	non--compact					 complete							 orbifold which											 is sufficiently
			 good at												 infinity. 
											 Let $D_i$ and											$D_i^\pm$ be the			 generalized												 Dirac operators
			 on the 
	 orbifold $X$ with											 coefficients in										 the Hermitian	 orbibundle
	 $E_i$ (with connection $\nabla^{E_i}$),  $i=1,2$									 Suppose that														 $D_1 = D_2$ in					 a neighborhood													 $\Omega$ of											 infinity. 
Assume that there exists a												 constant
	$k_0>0$ such that											$$
														 {\Cal R}\, \geq \,
														 k_0\, Id\, \hbox{
														 on } \Omega.
														 $$   
														 Then we know
														 from Theorem
														 3.2 that
														 ker$(D_i^\pm) <
														 + \infty$.
														 Thus the analytical index, of $D_i^\pm$, dimension of kernel minus dimension of cokernel,  is finite, i.e.,														 $ind_a(D_i^\pm)
														 < + \infty$,
														 $i=1,2$.

														 \proclaim{Definition	3.6} Let
														 $X$ be a
														 non--compact
														 complete												orbifold which
														 is sufficiently
														 good at
														 infinity. Let
														 ${ 
														 D}_i^+$,
														 ${ E}_i$,
														 $i=1,2$ be
														 generalized
														 Dirac
														 operators, with
														 coefficients in
														 the Hermitian														 orbibundles
														 $E_i$ (with connection $\nabla^{E_i}$), that
														 agree on a
														 neighborhood of
														 infinity
														 $\Omega$.
														 Then we define
														 the analytical
														 index, $ind_a(
														 { D}_1^+,
														 {D}_2^+)$,  of
														 the pair
														 $( { D}_1^+, {
														 D}_2^+)$  to be
														 $$
														 ind_a ( {														 D}_1^+, { 
														 D}_2^+) =
														 ind_a( {														 D}_2^+) -
														 ind_a( {														 D}_1^+),
														 $$
														 where $ ind_a(
														 { 
														 D}_i^+)$,
														 is the
														 analytical
														 index,
														 dimension of
														 kernel minus
														 dimension of
														 cokernel, of $D_i^+$,
														 $i=1,2$, see
														 Theorem
														 3.2. 
														 \endproclaim

								 \vskip 1em
								 \noindent {\bf 4. The Topological Index.}
								 \vskip 1em

								 We will define here the topological
								 index of a pair of generalized Dirac operators that
								 agree at infinity, thus generalizing a
								 construction of Gromov and Lawson,
								 [GL].
								 We start with a non--compact
					complete orbifold $X$ 
								 which is sufficiently
								 regular at infinity. 

								 Let $D_i$  and $D_i^\pm$ be the generalized 
								 Dirac operators on $X$ with coefficients
								 in
								 the Hermitian  orbibundle $E_i$ (with connection $\nabla^{E_i}$),  $i=1,2$.
								 Suppose that $D_1 = D_2$ in a
								 neighborhood of infinity. 

								 We assume here that the orbifold $X$ is
								 sufficiently good at infinity. Recall that this
								 means that, for any neighborhood of
								 infinity $\Omega$, we can chop off $X$
								 along an orbifold hypersurface,
								 $O_\Omega $, which is the boundary of a
								 neighborhood of infinity included in
								 $\Omega$. We also assume that the
								 orbifold structure of $X$ is of product
								 type in a neighborhood of $O_\Omega $.
								 Note that sufficiently good at infinity
								 implies sufficiently regular at
								 infinity.

								 Chop off $X$ along an hypersurface $O$,
								 with $D_1 = D_2$ on the neighborhood of
								 infinity of which $O$ is boundary. Glue
								 in an  orbifold along $O$, so that the
								 resulting space is a closed
								 orbifold $W$. (This can be
								 achieved, for example, by gluing
								 another copy of the chopped off
								 orbifold along $O$). We can assume that
								 $E_1=E_2$ on  a neighborhood of $O$. Extend $E_i$ to $W$, and call ${\breve E}_i$ this extension, $i=1,2$. We can assume that 
								 on the glued in part,  ${\breve E_1} =
								 {\breve E}_2$.
								 Let  ${\breve D}_i$ and
								 ${\breve D}^\pm_i$
								 be the generalized Dirac operators on
								 $W$ with coefficients in ${\breve 
								 E}_i$, $i=1,2$.
Define the topological index of the 							 pair $( { D}_1^+, { D}_2^+)$ by

								 \proclaim{Definition 4.1} Let $X$
								 be a non--compact complete
								 orbifold which is sufficiently good at
								 infinity. Let   $E_i$, ${D}_i$, ${ D}_i^+$, ${\breve E}_i$,   ${\breve D}_i^+$, $i=1,2$  be as above.
								 Then we define the topological index,
								 $ind_t( { D}_1^+, { D}_2^+)$,  of the pair $( { D}_1^+, { D}_2^+)$  to be
								 $$
 ind_t ( {D}_1^+, { D}_2^+) = ind_t(
								 {\breve D}_2^+  ) - ind_t( {\breve  D}_1^+),
								 $$
 where $ ind_t( {\breve D}_i^+)$,
 $i=1,2$ is the topological index of
 Kawasaki of the operator ${\breve D}_i^+$ on the closed orbifold $W$,
 $i=1,2$, [Kw2], [Kw3]. 
 \endproclaim

 We will now give an explicit local 
description of the above index by using
Kawasaki's formulas, see [Kw2], [Kw3], [Du], [V].
Let $\Cal U  =\{ ({\tilde U_i}, G_i )| i
\in { I\} }$, with ${\tilde U}_i/G_i =
U_i$ and projection $\pi_i: {\tilde
U}_i \to U_i$, be a standard orbifold
atlas of $W$.  Now ${\tilde U}_i^g $, the set of points fixed by $g$ in ${\tilde U}_i $, 
for any $g\in G_i$, admits the action
of the centralizer $Z_{G_i }(g)$ of $g$
in $G_i$. If $g$ and $g'$ are conjugate in 
$G_i$, then ${\tilde U}_i^g $ and
 ${\tilde U}_i^{g'} $  are diffemorphic
								 via some element $h$ in
								 $G_i$, with $g' = h g h^{-1}$. So we
								 can consider only one element for each
								 conjugacy class in $G_i$. For each
								 point $x\in W$, let
								 $(1), \dots, (h^{\rho_x}_x)$ be all the
								 conjugacy classes of the stabilizer
								 $G_x$ of $x$. Then we have a natural associated 	 orbifold,
[Kw3], 
$$
\hat{\Sigma}W=
\left\{ (y,
														(h^j_y))\, | \, y\in
														W, G_y \not=
														1,\, j=2, \dots,
														\rho_y
														\right\}.
														$$
														$\hat{\Sigma}W$
														is stratified by
														orbit types;
														define														$$
														\hat{\hat{\Sigma}}W=
														W \cup
														\hat{\Sigma}W.
														$$
														For short, we
														can rewrite 
														$\hat{\hat{\Sigma}}W$
														as
														$$
														\hat{\hat{\Sigma}}W
														= \coprod_{j=1}^q
														\hat{\Sigma}^jW,
														$$  
														where
														$\hat{\Sigma}^1W=W$,
														and so
														$\hat{\Sigma}^jW$
														is the 
														stratum
														corresponding to
														the $j$--th
														orbit type, $j=1,\dots, \ell$.
														In general, the
														action of
														$Z_{G_x}(h)$ on
														${\tilde
														U}_x^{h}$ is not
														effective. The
														order of the
														subgroup of
														$Z_{G_x}(h)$
														that acts
														trivially is
														called the
														multiplicity of
														$\hat{\Sigma}W$
														at $(x, h)$, and
														it is denoted by
														m ($x\in W$ and
														$h\in G_x$).
														Hence to each
														connected
														component of
														$\hat{\Sigma}W$
														we assign a
														certain constant
														multiplicity. 
														The following
														theorem, the
														index theorem
														for generalized
														Dirac
														operators on
														closed
														orbifolds, was
														proved by
														Kawasaki in [Kw1], 
		[Kw2]; see also [Du; Theorem 14.1], [BGV], [V].
														\proclaim{Theorem
														4.2 [Kw2], [Du]} Let
														$W$ be a closed 
														orbifold. Let
														$E$ be a Hermitian										orbibundle (with connection $\nabla^E$)  on
														$W$ and let
														${ D}_E^+$
														be a generalized
														Dirac operator
														with coefficients in
														$E$. Then the
														(topological and
														analytical)
														index of
														${ D}_E^+$
														on $W$ is given
														in terms of
														local traces by  
														$$
														\hbox{ind}({ D}_E^+) = \int_{\hat{\hat{\Sigma}}W } d\mu_{D_E^+}^\Sigma :=\sum_{j=0}^\ell \frac1{m_j} \int_{{\hat \Sigma}^j W} d\mu^{\Sigma, j}_{D_E^+},
	$$
														 where $
														 d\mu^{\Sigma,
														 j}_{D_E^+}
														 $ is the
														 density on
														 $j$-th stratum
														 of
														 $\hat{\hat{\Sigma}}W$
														 associated (via
														 parametrices)
														 to the operator ${  D}_E^+$, and
														 $m_j$ is the
														 corresponding
														 multiplicity
														 function,
														 $j=1,\dots,
														 \ell$. (For
														 simplicity's
														 sake, we
														 assumed all the
														 strata to be
														 connected; in
														 the general
														 case there will
														 be an
														 additional
														 summation over
														 the connected
														 components of
														 the strata.) 
														 \endproclaim
														 \proclaim{Remark
														 4.3 }
														 For 
														 another 
														 formulation of
														 Theorem 4.2
														 in terms of
														 orbifold
														 characteristic classes,
														 see [Kw3]. 
														 \endproclaim
														 We will now use
														 Theorem 4.2
														 to compute the
														 topological
														 index of a pair
														 of generalized
														 Dirac
														 operators on a
														 compact
														 orbifold that
														 is sufficiently
														 good at
														 infinity.
														 \proclaim{Theorem
														 4.4} Let $X$, $W$, ${D}_i$, 										 ${D}_i^+$, ${\breve D}_i^+$, ${E}_i$,
	 ${\breve E}_i$,
	 $i=1,2$ be as
	in Definition 4.1.  Let											 $Q_i$ be a	semi--local										 parametrix for	 ${ D}_i^+$,									 $i=1,2$, on $W$	such								 that $Q_1 =Q_2$	 in a									 neighborhood of
	 infinity; write
														 $$
														 D_i^+ Q_i =Id -														 T_i, \hbox{ and
														 } Q_i D_i^+
														 =Id - T'_i,
														 \,i=1,2,
														 $$
														 with $ T_i $
														 and $T_i'$ the
														 associated
														 semi--local
														 smoothing
														 operators. Then
														 the local trace
														 functions of
														 Kawasaki
														 associated to
														 $T_1$ and $T_2$
														 (and $T_1'$ and
														 $T_2'$ )
														 coincide at
														 infinity, and
														 the topological
														 index $ind_t (
														 { D}_1^+,
														 { 
														 D}_2^+)$ of the
														 pair $ (
														 { D}_1^+,
														 { D}_2^+)
														 $
														 is given by, 
														 $$
														 \hbox{ind}({\tilde D}_E^+)=
\int_{\hat{\hat{\Sigma}}X } \left(d\mu_{T'_2}^\Sigma- d\mu_{T_2}^\Sigma \right)- \int_{\hat{\hat{\Sigma}}X } \left( d\mu_{T'_1}^\Sigma- d\mu_{T_1}^\Sigma
\right)
														 $$ 
														 \endproclaim
														 \noindent {\it 
														 Proof.} This
														 theorem follows
														 directly from
														 Theorem
														 4.2.
														 Indeed, if we
														 denote by $\Omega$ the
neighborhood of infinity on where $D_1$ and $D_2$
coincide, we can cap off $X$ along an orbifold $O$ in
$\Omega$.
														 Consider a
														 parametrix
														 $ Q_0$ for
														 the extension
														 ${\breve
														 D}_1={\breve
														 D}_2$ on
														 $\Omega$. Now
														 splice  $Q_0$ onto
														 $Q_1$  and
														 $Q_2$, via a
														 smooth function
														 $f$ which is 0
														 outside
														 $\Omega$, 1
														 on a
														 neighborhood of
														 infinity, and
														 whose gradient
														 is bounded by 1
														 in norm. For
														 the rest of the
														 proof we
														 can proceed
														 exactly as in
														 [GL; Proof of
														 Proposition
														 4.6]. Note that we need to use Kawasaki's local trace formulas. $\qed$
														 \proclaim{Remark
														 4.5} The
														 topological
														 index of the
														 pair, $ind_t (
														 { D}_1^+,
														 {														 D}_2^+)$, is
														 independent of
														 the extension.
														 \endproclaim
														 \vskip 1em
														 \noindent {\bf
														 5. The
														 Analytical
														 Index:
														 Computations.}
														 \vskip 1em
														 We will now detail
														 some local
														 properties of
														 the  trace of a 
 generalized Dirac operator $D$, thus extending to
														 orbifolds some
														 results of
														 Gromov and
														 Lawson, see [GL].
														 Assume that
														 there exists a
														 constant
														 $k_0>0$ such
														 that
														 $$
														 {\Cal R}\,  \geq \, 
														 k_0\,  Id\, \hbox{
														 on } \Omega, 														 $$   
with $\Omega= X-K$ as before.
														  Because of Theorem 2.8, our assumption implies 
														$$ \Vert D\sigma \Vert_X \geq c \Vert \sigma \Vert_X, \, \forall \sigma \in {\Cal C}^{\infty}_c(X, {\Cal S}) \hbox{ such that } supp(\sigma) \cap K = \emptyset \tag 5.1 $$
												Our goal in this section is
														 to explicitly
														 compute the
														 index trace of 
														 $D$
														 in
														 terms of local
														 data, by
														 using
														 techniques of
														 Gromov and
											Lawson and
														 Anghel, see
														 [GL],
														 [An2].
														 This will
														 enable us to
														 prove, in
														 Section 6,
														 the relative
														 index theorem by localization.
														 Also, for
														 simplicity's
														 sake, in this
														 section we will
														 take
														 $E=\Bbb C$. So
														 $D$ will be
														 defined on
														 sections of the
														 $Spin^c$
														 bundle $\Cal
														 S$.
														 Start with a
														 parametrix
														 $Q_0$ of $D$ on
														 $X$. Then  
														 $$
														 D Q_0 = Id -
														 R, \hbox{ and
														 } Q_0 D  =Id -
														 R', \hbox{ on
														 $X$, }
														 $$
														 where $R$ and
														 $R'$ are not
														 necessarily
														 trace class (as
														 instead happens
														 in the closed
														 orbifold case).
														 We will now
														 replace $Q_0$,
														 outside a
														 compact set $K_0$, 
														 with a Green
														 operator
														 associated to a
														 suitable
														 extension of
														 $D$. Specifically, set $\Omega_0= X-K_0$, with $\Omega_0$ a domain with smooth boundary
included in $\Omega$. 
														 Let $D_{\Omega_0}$
														 be the graph
														 closure of the
														 restriction of
														 $D$ to ${\Cal
														 C}^{\infty}_c(\Omega_0,
														 {\Cal S})$ in
														 ${\Cal
														 L}^{2}(X, {\Cal
														 S})$.  
														 Let $P_{\Omega_0}$
														 denote the
														 orthogonal
														 projection from ${\Cal
														 L}^{2}(X, {\Cal
														 S})$ to 
	$$
														 H_{\Omega_0}
														 =\left\{ \sigma
														 \in {\Cal
														 L}^{2}(\Omega_0, 
														 {\Cal S}) |
														 D|_{\Omega_0}
														 (\sigma|_{\Omega_0})=0
														 \hbox{
														 distributionally
														 } \right\} 
														 $$
														 Note that, by
														 (5.1),
														 $H_{\Omega_0}$ is a
														 closed subspace
														 of $ {\Cal
														 L}^{2}(X, {\Cal
														 S})$. 
														 \proclaim{Theorem
														 5.1} Let
														 $X$ be a
														 non--compact
														 complete
														 orbifold which
														 is sufficiently
														 good at
														 infinity. Let
														 $D$, ${\Cal
														 S}$,
														 $D_{\Omega_0}$,
														 $P_{\Omega_0}$, be
														 as above.  Then
														 the following
														 two equations 
														 $$
														 D_{\Omega_0}
														 G_{\Omega_0} = Id
														 -P_{\Omega_0}\hbox{
														 on }{\Cal
														 L}^{2}(X, {\Cal
														 S})\tag 5.2														 $$
														 $$
														 G_{\Omega_0}
														 D_{\Omega_0} = Id
														 \hbox{ on the
														 minimal domain
														 of }D_{\Omega_0}
														 \tag 5.3$$
														 define a
														 bounded
														 operator
														 $G_{\Omega_0}$ :
														 ${\Cal
														 L}^{2}(X, {\Cal
														 S}) \to {\Cal
														 L}^{2}(X, {\Cal
														 S}) $, the
														 Green
														 operator of
														 $D_{\Omega_0}$.   
														 \endproclaim
														 \noindent{\it
														 Proof.}
														 Equation
														 (5.1)
														 shows that
														 $D_{\Omega_0}$ is
														 1:1 and has
														 closed range,
														 call it
														 ran$(D_{\Omega_0})$.
														 Hence  ${\Cal
														 L}^{2}(X, {\Cal
														 S})$ decomposes
														 orthogonally as
														 $$
														 {\Cal L}^{2}(X,
														 {\Cal S}) =
														 H_{\Omega_0} \oplus
														 \hbox{
														 ran}(D_{\Omega_0}).
														 $$
														 Thus $G_{\Omega_0}$
														 can be defined
														 to be zero on $
														 H_{\Omega_0}$ and $
														 D_{\Omega_0}^{-1}$
														 on
														 ran$(D_{\Omega_0})$.
														 $G_{\Omega_0}$ is
														 bounded because
														 of (5.1).
														 $\qed$
														 Since
														 $H_{\Omega_0}$ is a
														 closed subspace
														 of ${\Cal
														 L}^{2}(X, {\Cal
														 S})$, we can define
														 the orthogonal
														 projection on
														 $H_{\Omega_0}$
														 Bergman kernel
														 operator
														 $$
														 {\Cal P}_H
														 (x,y) =
														 \sum_m
														 \sigma_m (x)
														 \otimes
														 \sigma_m^*(y),
														 \, \forall
														 x,y,\in {\Omega_0},
														 $$
														 where $\{
														 \sigma_m\}$,
														 $m\in \Bbb N$,
														 is  an
														 orthonormal
														 basis
														 of $H_{\Omega_0}$.
														 \proclaim{Theorem
														 5.2} Let
														 $X$ be a
														 non--compact
														 complete
														 orbifold which
														 is sufficiently
														 good at
														 infinity. Let
														 $D$, ${\Cal
														 S}$,
														 $D_{\Omega_0}$,
														 $P_{\Omega_0}$, be
														 as above.  Then
														 the following
														 Bergman
														 kernel 
														 $$
														 {\Cal P}_H
														 (x,y) =
														 \sum_m
														 \sigma_m (x)
														 \otimes
														 \sigma_m^*(y),
														 \, \forall
														 x,y,\in {\Omega_0}
														 $$
														 of the
														 projection
														 operator
														 $P_{\Omega_0}$
														 converges (on
														 lifts of local
														 charts; see the
														 proof for
														 details)
														 uniformly in
														 the norm ${\Cal
														 C}^k$ on compact subsets, $k\in \Bbb N$.
														  \endproclaim
														  \noindent{\it
														  Proof.} (In
														  this proof
														  $\tilde {}$
														  means lifted to 
														  the local
														  orbifold
														  charts.)
														  Because $X$ is
														  sufficiently
														  good at
														  infinity, we
														  can assume
														  that there
														  exists a
														  smooth
														  exhaustion
														  function $F$
														  on $X$, such
														  that
														  ${\Omega_0}=
														  F(t_0)=\{
														  x\in X | F(x)
														  > t_0 \}$.  To
														  show the
														  convergence of
														  the Bergman
														  kernel ${\Cal
														  P}_H$ we
														  proceed in the
														  following way.
														  First of all,
														  set ${\Cal
														  B}_N=\sum_{m=1}^N \sigma_m(x) \otimes \sigma_m^*(y)$.
														  Obviously 
														  $B_N$ is a
														  finite rank
														  operator in
														  ${\Cal
														  L}^{2}(X,
														  {\Cal S})$.
														  Now, for any
														  $\phi \in
														  {\Cal
														  L}^{2}(X,
														  {\Cal S})$,
														  $B_N(\phi) \to
														  B(\phi)$. In
														  particular 
														  $B_N(\phi) \to
														  B(\phi)$ in
														  $\Cal D'$ (as
														  distributions).
														  By using local
														  charts,  $B_N
														  \to B$ weakly
														  in ${\Cal
														  L}_{G_U}
														  ({\Cal
														  D}_{\tilde y},
														  {\Cal D}'_{\tilde
														  x})$ for any
														  two points
														  $\tilde x$,
														  $\tilde y$ of
														  an orbifold
														  chart 
														  $(\tilde U,
														  G_U)$
														  projecting to
														  $x$ and $y$
														  respectively.
														  (By using a
														  partition of
														  unity, the
														  convergence of
														  distributions
														  is a local
														  property; use
														  a standard
														  orbifold
														  atlas.)  Also,
														  for
														  simplicity's
														  sake, we have
														  taken
														  $x$
														  and $y$ to  be
														  in
														  the same
														  orbifold chart.
														  Now,
														  $P_{\Omega_0}$,
														  which in this
														  proof we will
														  call $P$, is a
														  bounded
														  operator on
														  ${\Cal
														  L}^{2}(X,
														  {\Cal S})$,
														  and so is also
														  continuous as
														  an operator
														  $\Cal D \to
														  \Cal D'$.
														  Hence $P$ has
														  a
														  $G_U$--invariant
														  distributional
														  kernel
														  $p(\tilde x,
														  \tilde y)$ in
														  the sense of
														  Schwartz, c.f.
														  [Schw1],
														  [Scw2],
	 [At1]. By the Schwartz
											  kernel									  theorem,   $B_N \to {\Cal P}_H$
	 locally as a												  $G_U$--invariant
	  distribution											  on $\tilde U											  \times \tilde											  U$. (This											  follows as in the manifold case, see [At1; pg.  51].) So 
$$
{\Cal P}_H(\tilde x,
\tilde y)=
														  \sum_{m=1}^{+ \infty}														  \sigma_m
														  (\tilde x)
														  \otimes
														  \sigma_m^*(\tilde
														  y)
														  $$
														  $G_U$--invariantly,
														  on local
														  charts.
														  Now note that
														  ${\Cal P}_H$ and $B_N$
														  satisfy the
														  elliptic
														  equation
														  ${\tilde
														  D}^2=0$.
														  Therefore by
														  applying Lemma
														  3.1														  locally, we
														  have the
														  required
														  uniform
														  convergence.
														  $\qed$

														  The kernel ${\Cal P}_H$ has the strong
														  finiteness
														  property proved
														  below.
														  Let $F: X\to
														  {\Bbb R}^+$ be
														  a smooth
														  exhaustion
														  function as
														  above.
 In particular
														  we assume that
														  ${\Omega_0}=
														  F^{-1} (t_0, +
														  \infty)$, $K=
														  F^{-1} [0,
														  t_0]$. Let 
														  $
														  X(t)= \{  x\in
														  X|  F(x) >t
														  \}
														  $;   
														  then ${\Omega_0} =
														  X(t_0)$. We 														  have, c.f. [GL;  Lemma
														  4.20]
														  for the
														  manifold
														  case,
														  \proclaim{Theorem
														  5.3} Let
														  $X$ be a
														  non--compact
														  complete
														  orbifold which
														  is
														  sufficiently
														  good at
														  infinity. Let
														  $D$, ${\Cal
														  S}$,
														  $D_{\Omega_0}$,
														  $P_{\Omega_0}$,
														  ${\Cal P}_H
														  (x,x)$ be as
														  above.  
														  Then for any
														  $t>t_0$, 
														  $$
														  \int_{X(t)}
														  {\Cal P}_H
														  (x,x) < +
														  \infty
														  $$
														  \endproclaim
														  \noindent{\it
														  Proof.} Fixed
														  $t>t_0$ ($t$ near
														  $t_0$), choose
														  $s$ so that
														  $t_0<s<t$, and
														  consider the
														  compact
														  "annulus"
														  $A=$closure$[X(s)-X(t)]$. Since  ${\Cal P}_H (x,x)$ converges uniformly in the ${\Cal C}^k$,
														  $k\in \Bbb N$
														  norm, on compact subsets,
														  we can assume
														  that it
														  converges in
														  the  Sobolev
														  norm on $A$,
														  that is, 
														  $$\sum_m
														  \Vert \sigma_m
														  \Vert_{1,A} <
														  + \infty. \tag
														  5.4$$ 
														  We now claim
														  that there
														  exists a
														  constant $c$
														  so that, for
														  each $\sigma
														  \in H_{\Omega_0}$,
														  $$
														  \Vert
														  \sigma_m\Vert_{1,
														  X(t)}^2 \leq
\,
														  c \,
														  \Vert \sigma
														  \Vert^2_{1,A}.
														  \tag 5.5														  $$      
														  To do so,
														  choose a
														  cut-off
														  function $f\in
														  {\Cal
														  C}^{\infty}({\Omega_0})$
														  such that
														  $0\leq f \leq
														  1$, $f=1$ on
														  $X(t)$, and
														  $f=0$ on
														  ${\Omega_0}
														  -X(s)$.
														  Clearly there
														  exists a
														  $c_0\in {\Bbb
														  R}^+$ such
														  that $\Vert
														  \nabla f
														  \Vert_X <
														  c_0$.
														  Applying the
														  local identity
														  $$
														  \nabla (f
														  \sigma )=
														  \nabla (f
														  )\sigma +
														  f\nabla \sigma
														  $$  
														  and $D^2=
														  \nabla^*
														  \nabla + \Cal
														  R$,  we have,
														  for every
														  $\sigma \in
														  H_{\Omega_0}$,
														  $$
														  0=(D^2_{\Omega_0}
														  \sigma, f^2
														  \sigma) =
														  ((\nabla^*
														  \nabla+ \Cal
														  R)\sigma, f^2
														  \sigma)
														  $$
														  $$
														  0=(\nabla^*
														  \nabla\sigma,
														  f^2 \sigma)
														  +(\Cal R
														  \sigma, f^2
														  \sigma)
														  $$
														  $$
														  0=(
														  \nabla\sigma,
														  \nabla(f^2 \sigma))
														  +(\Cal R
														  \sigma, f^2
														  \sigma)
														  $$
														  $$
														  0=(\nabla\sigma,
														  2 f (\nabla f)
														  \sigma) +
														  (\nabla\sigma,
														  f^2 (\nabla
														  \sigma))
														  +(\Cal R
														  \sigma, f^2
														  \sigma)
														  $$
														  $$
														  0=2 (f
														  \nabla\sigma,
														  (\nabla f)
														  \sigma) + ( f
														  \nabla\sigma,
														  f (\nabla
														  \sigma))
														  +(\Cal R
														  \sigma, f^2
														  \sigma). 
														  $$
														  Since $<\Cal R
														  \alpha,
														  \alpha> \, \geq \, 
														  k_0 \,  <\alpha,
														  \alpha>$ in
														  ${\Omega_0}$, we
														  have
														  $$
														  \Vert \nabla
														  \sigma
														  \Vert^2_{X(t)}
														  + k_0  \Vert
														  \sigma
														  \Vert^2_{X(t)}
														  $$
														  $$
														  \leq 2 |(f
														  \nabla\sigma,
														  (\nabla f)
														  \sigma)| \leq
														  2c_0 \Vert
														  \nabla \sigma
														  \Vert_A
														  \Vert \sigma
														  \Vert_A \leq
														  c_0 (\Vert
														  \nabla \sigma
														  \Vert^2_A
														  +\Vert \sigma
														  \Vert^2_A).
														  $$
														  Then
														  $$
														  \Vert \sigma
														  \Vert^2_{1,
														  X(t)} \leq c_0
														  (1+\frac1{k_0})
														  \Vert \sigma
														  \Vert^2_{1, A}.
														  $$
														  We thus proved
														  our claim
														  (5.5).
														  Now the Lemma
														  is proved by
														  combining
														  (5.4) and
														  (5.5). $\qed$
\vskip 1em
\noindent {\bf 6. The Relative Index Theorem.}
\vskip 1em
We will prove here that the analytical index of a pair of generalized Dirac
operators that agree at infinity is equal to its analytical index, thus
generalizing to orbifolds the main theorem of Gromov and Lawson in
[GL].
Let $X$ be a non--compact complete orbifold which is
sufficiently good at infinity. Assume that there exists a constant
$k_0>0$ such that
$$
{\Cal R} \geq k_0 Id\, \hbox{ on } \Omega, \quad K=X-\Omega,
$$   
with $\Omega$ a domain with smooth boundary, and $\Cal R$ as in Proposition 2.7.
Let $D_i$  and $D_i^\pm$ be the
generalized Dirac operators on $X$ with
coefficients in the Hermitian orbibundle $E_i$ (with connection $\nabla^E_i$),  $i=1,2$. Suppose that $D_1
= D_2$ in a neighborhood $\Omega$ of infinity.
Let
$G_i$, $i=1,2$ be the Green operator associated to $D_i^+$, $i=1,2$, see
Theorem 3.5. In
particular, $D_i^+$ and $G_i$, satisfy the following equations,
$$
D_i^+ G_i = Id - P_i^-, \hbox{ and } G_i D_i^+ = Id-P_i^+,\quad i=1,2,
$$
where $P_i^\pm : {\Cal L}^{2}(X, {\Cal S}) \to  {\Cal L}^{2}(X, {\Cal S})$
is the projection onto the finite dimensional space ker$ D_i^\pm$.
On $\Omega$, $D_1 = D_2$. Restrict $G_i$ to
$\Omega$
by defining ${\hat G}_i = \chi{G}_i \chi
$, $i=1,2$, where  $\chi $ is the characteristic function of $\Omega$.
The difference ${\hat G}_2  - {\hat G}_1$ satisfies
$$
D^+ ({\hat G}_2  - {\hat G}_1) = \hat{P}_1^- - \hat{P}_2^-,
$$
where $D_1 =D_2 =D$ on $\Omega$, and $\hat{P}_i^\pm = \chi P_i^\pm \chi $, $i =
1,2$,
is a finite range operator. This implies, as in [GL], that the range of
$({\hat G}_2^+  - {\hat G}_1^+) $  is nearly contained in the kernel
$H_\Omega$
of $D^+$ on $\Omega$, notation as in Proposition 2.7. In other words, 
let $V$ be ker$({\hat P}_2^+  - {\hat P}_1^+) $. Then $V$ is a subspace
of finite codimension in ${\Cal L}^2 (X, {\Cal S}^+)$, and
$$
({\hat G}_2  - {\hat G}_1) V \subseteq \hbox{ ker}(D^+)
$$  
\proclaim{Theorem 6.1} (.f. [GL; Lemma 4.28] for the manifold
case.)Let $X$,  $ D_i$,  $ G_i$, $ {\hat G}_i$,
$i=1,2$, 
be as above. Then the local trace density function of  $({\hat G}_2  -
{\hat G}_1)$ is integrable at infinity. 
\endproclaim
\noindent{\it Proof.} Recall that  $\Omega= X-K$, with $K$ compact. Choose 
$\Omega'= X-L$, with $L$ compact such that $K\subseteq $int$(L)$. By
Theorem 5.3, if $P_\Omega$ denotes the orthogonal projection onto the
kernel $H_\Omega$ of $D^+$ on $\Omega$, and ${\Cal P}_\Omega$ is its
associated trace density, then
$$
\int_{\Omega'} {\Cal P}_\Omega (x) \, dv \, < \, + \infty.
$$ 
Set $Z=({\hat G}_2  - {\hat G}_1) $. Then
$$
\hbox{range}(Z) \subseteq H_\Omega +F,
$$
where $F$ is a finite--dimensional subspace of ${\Cal L}^2 (X, {\Cal
S}^-)$,
by the discussion preceding the statement of Theorem 6.1. Let
$\{\sigma_m\}$, $m\in \Bbb N$, be an orthonormal basis of   $ H_\Omega
+F$, 
such that $\{\sigma_m\}$, $m= N, N+1, \dots $, is an orthonormal basis
of 
$H_\Omega$. Then the Schwartzian kernel of $Z$ can be written as 
$$
K^Z(x,y) = \sum_m \sigma_m(x) \otimes (Z^*\sigma_m)(y), \quad \forall x,y, \in \Omega,
$$  
where $Z^*$ denotes the adjoint of $Z$. The local trace function ${\Cal
P}_Z$ of $Z$ satisfies 
$$
|{\Cal P}_Z (x)| \leq \sum_m | <\sigma_m(x) , (Z^*\sigma_m)(x)>| , \quad
\forall x \in \Omega.
$$
Let $Z'= \chi'  Z\chi' $, with $\chi'$ the characteristic  function of
$\Omega'$, denote the restriction of $Z$ to $\Omega'$; 
note that
$\Vert Z' \Vert \leq \Vert Z \Vert$. Then
$$
\int_{\Omega'}|{\Cal P}_Z (x)| \leq \sum_m  \int_{\Omega'} |<Z
\sigma_m(x) ,\sigma_m(x)>|\, dx 
$$
$$
\leq \, \sum_m  \Vert  Z \sigma_m \Vert_{\Omega'}\, \Vert
\sigma_m\Vert_{\Omega'}  
$$
$$
\leq \,  \sum_m  \Vert  Z' \sigma_m \Vert_{\Omega'}\, \Vert
\sigma_m\Vert_{\Omega'}  
$$
$$
\leq \, \Vert Z\Vert \sum_m  \Vert (\sigma_m)\Vert_{\Omega'}^2  
$$
$$
\leq \, \Vert Z\Vert \left( \sum_m^{N-1}  \Vert
(\sigma_m)\Vert_{\Omega'}^2 +   
\int_{\Omega'}|{\Cal P}_Z | \right)\, < \, + \infty \quad \qed
$$
We are now in a position to state and prove our relative index theorem.
\proclaim{Theorem 6.2} Let $X$ be a
non--compact
complete orbifold which is sufficiently good at infinity. Let $D_i$,
$D_i^\pm$ be the
generalized Dirac operators on $X$ with coefficients in the
Hermitian orbibundle $E_i$ (with connection $\nabla^E_i$),  $i=1,2$. Suppose that $D_1 = D_2$ in a
neighborhood $\Omega$ of infinity with smooth boundary. 
Assume that there exists a constant $k_0>0$ such that
$$
{\Cal R}\, \geq \,  k_0 \, Id\, \hbox{ on } \Omega.
$$   
where ${\Cal R}$ is given as in Section 5.
Then the topological index (as
in Definition 4.1) and the analytical index (as in Definition
3.6) of the pair $( {D}_1^+, { D}_2^+)$ coincide, that is

$$
ind_t ( {D}_1^+, { D}_2^+) = ind_a ( {D}_1^+, { D}_2^+) $$
\endproclaim 
The proof of Theorem 6.2 will occupy the rest of this section.

\noindent{\it Proof.(c.f. [GL; Proof of Theorem 4.18]) }
We will now construct parametrices $Q_j$ for $D^+_j$, $j=1,2$, by 
cutting off the Green' operators $G_i$, $i= 1,2$ of Theorem 3.5 in a neighborhood of
the diagonal. (Notation as at the beginning of this section.)  
Choose $\psi\in {\Cal C}^{\infty} (X\times X)$ with support in a small
neighborhood of the diagonal, so that $0\leq \psi \leq 1$ and $\psi = 1$
near the diagonal.  Note that $\psi$ can in particular be lifted, one
variable at the time, to a $G_U$--invariant function in $\tilde \psi$ in
${\Cal C}^{\infty}({\tilde U} \times {\tilde U})$, for any local
orbifold chart $(\tilde U, G_U)$. Let ${ R}_i$, $i=1,2$, be the
operator whose Schwartzian kernel is locally given by
$$
K^i (\tilde x, \tilde y)= \tilde \psi \, K^{G_i} (\tilde x, \tilde y)
$$ 
The operator ${ R}_i$, $i=1,2$, is a semi--local parametrix for
$D_i^+$
with
$$
D_i^+ R_i = Id - S_i^-, \quad R_i D_i^+ = Id - S_i^+,\quad i=1,2.
$$
As $R_i = G_i $ near the diagonal, the local traces associated to
$S^\pm_i$, $i=1,2$, satisfy 
$$
t^{ S^+_i} = t^{ i^+}, \quad t^{ S^-_i} = t^{ i^-},
$$
where $ t^{ i^\pm}$ is the trace associated to   $P_i : {\Cal L}^{2}(X,
{\Cal S}) \to  {\Cal L}^{2}(X, {\Cal S})$,
the projection onto the finite dimensional space ker$ D_i^+ $ and ker$
D_i^- $=
coker$ D_i^+ $, $i=1,2$.
From Theorem 4.4, to compute the topological index we need a pair of
semi--local parametrics that agree in a neighborhood of infinity. We will
do this by splicing $R_2$ onto $R_1$ in $\Omega$ as follows. 
Let $\{ b_k \}$ $k\in \Bbb N$, be a sequence of functions as in
Proposition 2.4. We assume that, for $k$ sufficiently large, $b_k=1$
on $K$. Moreover, since supp$(b_k)\subseteq \overline{B}_{2k}$, we can
assume that $b_k=0$ in a neighborhood of infinity. We claim that each
$b_k$ can be approximated by a smooth function $f_k$ such that  $f_k=1$
on $K$, supp$(f_k) \subseteq \overline{B}_{2k}$. This claim is proven in
the following way. Supp$(b_k)$
can be covered by the union of finitely many orbifold charts, say
${U}_j$ with 
 ${\tilde U}_j/G_j=U_j$, $j=1, \dots, \ell$. If ${\tilde
b}_{k,s}$ is the lift of a smooth approximation of $b_k$ on $U_s$, 
obtained $G_s$--invariantly, then
${\tilde f}_k =\sum_s {\tilde \eta}_s {\tilde b}_{k,s}$  is a smooth
approximation of $b_k$, $k\in \Bbb N$.      
Now define a sequence of semi--local parametrices $R_{2,k}
$ for $D_2^+$ by
$$
R_{2,k} = f_k R_2 + (1-f_k) R_1
$$ 
(recall that on $\Omega$,  $D_1^+=D_2^+$.)
Moreover,
$$
D_2^+ R_{2,k} =Id -S_{2,k}^-,\quad R_{2,k} D_2^+ =Id -S_{2,k}^+,
$$
with 
$$
S_{2,k}^- =f_k S_{2}^- + (1-f_k) S_1^- + \nabla(f_k) (R_2-R_1)
$$
$$
S_{2,k}^+ =f_k S_{2}^+ + (1-f_k) S_1^+ \tag 6.1$$
where $S_i^\pm$, $i=1,2$,  satisfy
$$
D_1^+ R_1 = Id -S_1^-, \quad R_1 D_1^+ = Id -S_1^+,
$$
$$
D_2^+ R_2 = Id -S_2^-, \quad R_2 D_2^+ = Id -S_2^+.
$$
Now by Theorem 4.4 applied to the pair of parametrices 
$R_{2,k}, R_1$, we get,
$$
ind_t ( {D}_1^+, { D}_2^+) = \int_{\hat{\hat{\Sigma}}X} \left( d\mu^{ S_{2,k}^+
} - d\mu^{ S_{2,k}^- } \right) 
-\int_{\hat{\hat{\Sigma}}X} \left( d\mu^{ S_{1}^+ } - d\mu^{ S_{1}^- }\right).  
$$
By (6.1) this latter expression is equal to
$$
=\int_{\hat{\hat{\Sigma}}X} f_k \left( d\mu^{ S_{2}^+ }-d\mu^{ S_{2}^-
} - d\mu^{ S_{1}^+ } +d\mu^{ S_{1}^- }\right) -
\int_{\hat{\hat{\Sigma}}\Omega}
d\mu^{\nabla f_k(R_2 .R_1)}
$$
$$
=\int_{\hat{\hat{\Sigma}}X} f_k \left( d\mu^{ P_{2}^+ } -d\mu^{ P_{2}^-
} - d\mu^{ P_{1}^+ } +d\mu^{ P_{1}^- }\right) -
\int_{\hat{\hat{\Sigma}}\Omega}
 d\mu^{\nabla f_k  (R_2 .R_1)},
$$
where $P_i^\pm$, is the projection onto the kernel of $D_i^\pm$, $i =1,2$.
Now the Schwartz kernel of $Z_k = \nabla f_k (R_2-R_1)$ on a local
chart is
$$
K^{Z_k} ( \tilde x, \tilde y) = \nabla \tilde{f}_k K^{R_2-R_1} ( \tilde
x, \tilde y).
$$
Near the diagonal, $R_i =G_i$, $i=1,2$, and so
$$
| d\mu ^{Z_k} (x)| \leq \Vert  \nabla f_k \Vert \,
| d\mu ^{Z} (x)|,
$$
where $Z= {\hat G}_2 -{\hat G}_1$ as above.
By integrating on ${\hat{\hat{\Sigma}}\Omega} $,  since supp$(\nabla
f_k)\subseteq \Omega$, we have
$$
\int_{\hat{\hat{\Sigma}}\Omega} d\mu^{\nabla f_k(R_2 -R_1)}=
\int_{\hat{\hat{\Sigma}}\Omega'} d\mu^{\nabla f_k(R_2 -R_1)}
$$ 
$$
\leq \int_{\hat{\hat{\Sigma}}\Omega'} \Vert \nabla f_k \Vert \,
| K^{R_2 -R_1}|
= \int_{\hat{\hat{\Sigma}}\Omega'} \Vert \nabla f_k\Vert \,
 | K^{Z}|.
$$
Since $\int_{\hat{\hat{\Sigma}}\Omega'} | K^{Z}| < + \infty$ by Lemma
5.3, our result follows if we show  that 
$$
\hbox{sup}_\Omega \Vert \nabla f_k \Vert \to 0 \hbox{ as } k\to + \infty.
$$
To show this latter claim, we will show that $sup_\Omega\Vert \nabla f_k\Vert \leq
\frac{M}k$ for $k \in \Bbb N$. For,  the compact
subset supp$(b_k)$ can be covered by finitely many orbifold charts. On each of
these charts, we can assume that we have an approximation with
sup$\Vert \nabla f_k \Vert \leq \frac{M}k $, $k\in \Bbb N$,  because
we first approximate on its  lift,  and then average over the chart
group.
Each function (and each  function's gradient) in the average will have the same sup.
Moreover, on the lift of each chart we can assume that  the gradient of the
distance function is bounded by  $\frac{M}k $ except on a (singular) set
of zero measure and codimension at least 2. Hence by possibly performing
a cutting and limiting procedure around the singular locus, a smooth
approximation with the required bound can be found. $\qed$  
\vskip 1em
\noindent {\bf References}
\itemitem {[A]} S. Agmon, Lectures on elliptic boundary value problems, 
D. Van Nostrand Co., Inc., Princeton, N.J.-Toronto-London, 1965. 
\itemitem {[An1]} N. Anghel, Extrinsic upper bounds for eigenvalues of
Dirac-type operators,  Proc. Amer. Math. Soc.  117  (1993),  501--509. 
\itemitem {[An2]} N. Anghel, An abstract index theorem on noncompact
Riemannian manifolds,  Houston J. Math.  19  (1993), 223--237. 
	   \itemitem {[At1]} M.F. Atiyah, Elliptic operators, discrete
	   groups and von Neumann algebras, Colloque "Analyse et Topologie"
	   en l'honneur de Henri Cartan (Orsay, 1974), pp. 43--72.
	   Asterisque, No. 32-33,
	   Soc. Math. France, Paris, 1976.
		   \itemitem {[BGV]} N. Berline, E. Getzler and M. Vergne, Heat
			   kernels and
				   Dirac operators, Grundleheren der Mathematical
				   Wissenshaften 298,
					   Springer—Verlag,  Berlin 1992.
					   \itemitem {[B]} J. Borzellino, Orbifolds of
					   maximal diameter,
					   Indiana Univ. Math. J. 42 (1993), 37--53.
					   \itemitem {[C]} G. Chen, Calculus on orbifolds,
					   Sichuan Daxue Xuebao  41  (2004),					   931--939. 
								\itemitem {[Ch]} Y.-C. Chiang, Harmonic
								maps of $V$-manifolds,           Ann.
								Global Anal. Geom.  8  (1990), 315-344.
								\itemitem {[Du]} J. J. Duistermaat, The
								heat kernel Lefschetz fixed point
								formula for the $Spin^c$ Dirac operator,
								Progress in Nonlinear Differential
								Equations and their Applications, 18.
								Birkh\"auser Boston, Inc., Boston, MA,
								1996. 
									   \itemitem {[Fa1]} C. Farsi,
									   $K$-theoretical index theorems
									   for
										   orbifolds,
											   Quat. J. Math. 43 (92),
											   183--200.
												   \itemitem {[Fa2]}
												   $\underline {\hskip
												   1.0in}$   Orbifold
												   spectral
													   theory,
														   Rocky Mtn. J. Math.
														   31 (2001), 
														   215--235.
														\itemitem {[Fa3]}$\underline {\hskip 1.0in}$ Orbifold $\eta$-invariants,
to appear.

http://www.iumj.indiana.edu/IUMJ/forthcoming.php
															   \itemitem
															   {[Fa4]}
															   $\underline
															   {\hskip
															   1.0in}$ Dirac operators on non--compact orbifolds, 
															   preprint
								 math.DG/0601451.
															   \itemitem
															   {[Fa5]}$\underline {\hskip 1.0in}$ Analysis on orbifolds, 
															   in
															   preparation.
															   \itemitem
															   {[GL]} 
															   M. Gromov
															   and  M.
															   Lawson, 
															   Positive
															   scalar
												   curvature
															   and the
															   Dirac
															   operator
															   on
															   complete
															   Riemannian
															   manifolds,
															   Inst.
															   Hautes
															   Etudes
															   Sci.
															   Publ.
															   Math. No.
															   58,
															   (1983),
															   83--196
															   (1984).
\itemitem {[Gn1]} M. Gaffney,  The harmonic operator for exterior
differential forms,  Proc. Nat. Acad. Sci. U. S. A.  37,  (1951).
48--50. 
\itemitem {[Gn2]} M. Gaffney, A special Stokes's theorem for complete
Riemannian manifolds,  Ann. of Math. (2)  60,  (1954), 140--145. 
\itemitem {[LR]} W. L\"uck and J. Rosenberg, Equivariant Euler  characteristics
and $K$-homology  Euler classes  for
																   proper
																   cocompact
																   $G$-manifolds, 
																   Geom.
																   Topol.
																   7
																   (2003),
																   569--613.
																   \itemitem
																   {[Ku]}
																   D.
																   Kucerovsky,
																   A
																   short
																   proof
																   of an
																   index
																   theorem,
																   Proc.
																   Amer.
																   Math.
																   Soc.
																   129
																   (2001),
																   3729--3736.
\itemitem {[Kw1]} T.  Kawasaki, The signature theorem for $V$--manifolds, Topology 17 (78), 75--83.  \itemitem {[Kw2]} 
$\underline {\hskip 1.0in}$ The Riemann Roch theorem for complex $V$--manifolds, Osaka J. Math. 16 (1979), 151--159.
\itemitem {[Kw3]} $\underline {\hskip 1.0in}$ The index of elliptic operators over $V$-manifolds, Nagoya Math. J. 84 (81), 135--157. 
\itemitem {[LM]} H. B. Lawson, Jr, and M.-L. Michelsohn, Spin Geometry, Princeton University Press, Princeton, New Jersey, 1989.
\itemitem {[N]} Y. Nakagawa, An isoperimetric inequality for orbifolds, Osaka J. Math. 30 (1993), 733--739. 							
\itemitem {[Sch1]} L. Schwartz, Theorie des distributions, Tomes I,
II. Act. Sci. Ind., no.  1091, 1122 = Publ. Inst. Math. Univ.
Strasbourg 9, 10  Hermann \& Cie., Paris, 1950,  1951.		
\itemitem {[Sch2]} L. Schwartz,
																   Laurent
																   Theorie
																   des
																   distributions à valeurs vectorielles. I.  Ann. Inst. Fourier, Grenoble  7  (1957), 
																   1--141. 
																   \itemitem
																   {[Sh]}
																   M. A.
																   Shubin,
																   Pseudodifferential operators and spectral theory. Translated from the 1978 Russian
																   original
																   by
																   Stig
																   I.
																   Andersson.
																   Second
																   edition,
																   Springer-Verlag, Berlin, 2001.
																   \itemitem
																   {[Si]}
																   J.
																   Simons,
																   Minimal
																   varieties
																   in
																   Riemannian
																   manifolds,
																   Ann.
																   of
																   Math.
																   (2)
																   88
																   (1968), 
																   62--105.
																   \itemitem
																   {[Stan]}
																   E.
																   Stanhope,
																   Spectral
																   bounds
																   on
																   orbifold
																   isotropy,  
																   Ann.
																   Global
																   Anal.
																   Geom.
																   27
																   (2005),
																   355--375.
																   \itemitem
																   {[Y]} S.
																   T. Yau, 
																   Some
																   function-theoretic properties of complete Riemannian manifold and their applications
																   to
																   geometry,
																   Indiana
																   Univ.
																   Math.
																   J.
																   25
																   (1976),
																   659--670.
											
\itemitem {[V]} M. Vergne, Equivariant index formulas for orbifolds, Duke Math J.  82 (1996), 637--652.  
																  \itemitem {[W]}																   J.
			  Wolf,
																   Essential
																   self-adjointness for the Dirac operator and its square,  Indiana Univ. Math. J.  22
																   (1972/73),
																   611--640.
																   \end